\documentclass[11pt]{article}

\newcommand{\N}{{\bf N}}

\newcommand{\RP}{{\bf R}P}
\newcommand{\C}{{\bf C}}
\newcommand{\X}{{\cal X}}
\newcommand{\M}{{\cal M}}

\newcommand{\inner}[2]{\left\langle#1,#2\right\rangle}%
\newtheorem{theorem}{Theorem}[section]

\newtheorem{lemma}[theorem]{Lemma}
\newtheorem{proposition}[theorem]{Proposition}
\newtheorem{definition}[theorem]{Definition}
\newtheorem{remark}[theorem]{Remark}%

\newcommand{\proof}[1]{\noindent{\it Proof#1:\  }}

\usepackage{amsfonts}
\title{ Positivity of Symplectic Area for Perturbed $J$-holomorphic
       Curves}
\author{Pawel Felcyn\\
        Department of Mathematical and Computer Sciences\\
        University of Wisconsin-Whitewater\\
        email: felcynp@mail.uww.edu}
\begin{document}

\maketitle
\begin{abstract}
\noindent
In this paper we will prove that for a compact, symplectic manifold
$(M, \omega)$  and for $\omega$-compatible almost-complex structure $J$
any properly perturbed $J$-holomorphic curve has a non-negative symplectic 
area.
This non-negative property provides us with a new obstruction to the
bubbling off phenomenon and thus allows us to redefine the Floer symplectic
homology.
In particular, in subsequent papers, we will prove the Arnold conjecture in
both degenerate and non-degenerate cases with integer coefficients for
general, symplectic manifolds.
\end{abstract}
\section{Introduction}
Let us denote by $(M,\omega)$ a compact, symplectic manifold. Here
$\omega$ is a closed, non-degenerate 2-form. If $H: \mathbb{R}
\times M \to \mathbb{R}$ denotes a time-dependent Hamiltonian
function then there is associated a time-dependent Hamiltonian
vector fields $X_{H_t}$ on the symplectic manifold $M$ defined by
the equality
   \begin{equation}
    \label{hamvector}
       \iota (X_{H_t})\omega = dH_t.
   \end{equation}
We shall assume that the Hamiltonian function $H$ is of period 1 in time.

The Arnold conjecture says that the number of 1-periodic solutions
of the Hamiltonian equation
   \begin{equation}
    \label{hameq}
     \dot{\gamma}(t)=X_{H_t}(\gamma (t))
   \end{equation}
is estimated from below by the number of critical points of a smooth
function defined over the manifold $M$. (See \cite[Appendix 9]{A}.)

In trying to prove the Arnold conjecture one is usually led to study
the action functional
  \[
     a_H(\gamma)= - \int_D u^* \omega +\int_0^1 H_t(\gamma(t)dt
  \]
defined on the space of smooth, contractible loops $\gamma:
 \mathbb{R} \to M$ with $\gamma(t)=\gamma(t+1)$. See, for example,
\cite{CZ,F5,H,HZ,S} for some special cases for both nondegenerate
and degenerate 1-periodic solutions of the equation (\ref{hameq}).

In trying to extend the variational methods to the general case one encounters
two difficulties:
  \begin{enumerate}
   \item
    The action functional $a_H$ is not uniquely defined on the space of
    smooth, contractible loops of $M$. Rather, it is well defined on the
    universal covering of the space of loops.
   \item
    There is bubbling off phenomenon which causes difficulties with
    compactification of appropriate moduli spaces.
  \end{enumerate}
In \cite{F2,F3,F4,HS,O1,O2,O3} these difficulties has been overcome in some
more general cases.

In our approach to the Arnold conjecture which is valid for all
compact symplectic manifolds we will restrict the action
functional to a subset of its one special `branch'. Thus the first
difficulty will be overcome in general. Fortunately, over the
restricted set the bubbling problem will also disappear. Thus, in
particular, we redefine the Floer symplectic homology to obtain a homology
theory very similar to the finite dimensional Morse homology.
Details of the construction of the new Floer homology and some of its 
consequences will be presented in subsequent papers.

In this paper we present the major tool in our construction of Floer
homology which is the positivity of symplectic area of properly perturbed
$J$-holomorphic curves. For a compact Riemann surface $(\Sigma,j)$ the
map $u:\Sigma \to M$ is said to satisfy a properly perturbed Cauchy-Riemann 
equation if
   \[
      du+J \circ du \circ j + P(u)=0.
   \]
(See Def \ref{percurve} for details.)
  \begin{theorem}
    Let $u:\Sigma \to M$ be a properly perturbed $J$-holomorphic curve.
    Then the symplectic area of $u$ is non-negative:
      \[
        \int_{\Sigma} u^* \omega \geq 0.
      \]
  \end{theorem}
 {\bf Acknowledgment} We are very thankful J.W. Robbin
 for finding mistakes in the previous version of these notes.
\section{Example}
Let $S^2$ be a Riemannian sphere $S^2=\C \cup \{\infty\}$ with the standard
Kaehler metric
   \[
     \frac{dz \otimes d{\overline{z}}}{(|z|^2+1)^2} ,
   \]
where $z=x+iy$ denotes a point of $S^2$. The induced standard symplectic
form is of the form
   \[
     \omega = \frac{dx \wedge dy}{(x^2 + y^2 +1)^2} ,
   \]
and the induced metric by the standard complex structure on $S^2$
is of the form
  \begin{equation}
  \label{metric}
    \inner{{\widehat{z}}_1}{{\widehat{z}}_2}=
    \frac{{\widehat{x}}_1{\widehat{x}}_2 +{\widehat{y}}_1{\widehat{y}}_2}{
          (|z|^2+1)^2},
  \end{equation}
where ${\widehat{z}}_j = {\widehat{x}}_j + i{\widehat{y}}_j$ denote
a tangent vectors at $z= x+iy$.

We will consider a Hamiltonian function $H: \mathbb{R} \times S^2
\to \mathbb{R}$ of the form
  \[
     H(s,z)= \psi(s) \frac{x^2+y^2-1}{x^2+y^2+1}.
  \]
The gradient vector field (depending on $s$) of the Hamiltonian
function $H$ with respect to the metric (\ref{metric}) is $\nabla
H(s,z) = 4\psi(s)z$. Thus we will obtain the following equation
for the properly perturbed holomorphic curve $u: \mathbb{R} \times S^1
\times S^2 \to S^2$ :
  \begin{equation}
  \label{simple}
   \frac{\partial u}{\partial s}+i\frac{\partial u}{\partial t}+4\psi(s)u=0.
  \end{equation}
Let us consider solutions of the equation (\ref{simple}) of the form
  \begin{equation}
   \label{solution}
    u(s,t)= \exp\left\{-4 \int_{-\infty}^s \psi(s')ds'\right\}
    \exp\{2 \pi k(s+it)\}.
  \end{equation}
Indeed, it is easy to verify that each function of the form (\ref{solution}),
for which the function $\int_{-\infty}^s \psi(s')ds'$ is smooth, satisfies
the perturbed Cauchy-Riemann equation (\ref{simple}).

Assume now that the function $\psi$ has a compact support. Under this
condition we claim that the symplectic area of solutions (\ref{solution})
is non-negative i.e.
   \[
     \int u^* \omega \geq 0.
   \]
Perhaps, the simplest way to see this is to homotop the curve $u(s,t)$
to a holomorphic one. The homotopy can be done by considering the equation
(\ref{simple}) with a parameter $\lambda$, $0 \leq \lambda \leq 1$ :
   \[
   \frac{\partial u}{\partial s}+i\frac{\partial u}{\partial t}+
   4 \lambda \psi(s)u=0.
   \]
Then for each $\lambda$ solutions of the above equation of the form
   \[
    u_{\lambda}(s,t)= \exp\left\{-4 \lambda \int_{-\infty}^s \psi(s')ds'
    \right\} \exp\{2 \pi k(s+it)\}
   \]
will provide us a homotopy between the original solution and the holomorphic
curve $u_0(s,t)=\exp{2 \pi k(s+it)}$. Now since the symplectic form
$\omega$ is closed and the symplectic area of the holomorphic curve is
non-negative by Stokes theorem we obtain
   \[
    \int u^* \omega = \int u^* \omega_0 \geq 0.
   \]
In fact, as we will see, the argument of deforming a properly perturbed
Cauchy-Riemann equation to a non-perturbed one works in general and thus
proving the positivity of symplectic area of perturbed $J$-holomorphic
curves.

Note, however, that if we let the function $\psi$ to be a nonzero
constant: $\psi(s)=\tau$ for all $s \in \mathbb{R}$, then as it
has been shown in \cite{HS} the function
  \[
    u_k(s,t)= \exp(4\tau s) \exp\{2\pi k (s+it)\}
  \]
is a solution of the equation (\ref{simple}) whenever $\pi k+2\tau>0$.
Moreover, for the symplectic area we have
  \[
    \int u_k^* \omega = \pi k.
  \]
Thus, in particular, we obtain solutions $u_k$ with a negative
symplectic area if $k$ is negative. 

The key point why the function $\psi$ with compact support produces only
solutions with non-negative symplectic area is that the Hamiltonian
part of the equation (\ref{simple}) represents a derivative of a global
function on $S^2$. In the case of the function $\psi$ being a nonzero 
constant this is not so.

\section{ Perturbed $J$-holomorphic Curves}\label{sec2}
Let $(M,\omega)$ denote a $2n$-dimensional compact symplectic manifold
 with the symplectic form $\omega$ and let $(\Sigma, j)$ denotes a closed
connected Riemman surface with a complex structure $j$ and with a fixed
Kaehler metric.
Let $J$ denote a smooth family of $\omega$-compatible almost-complex structure
 on $M$  depending on the parameter $z \in \Sigma $.
We will denote the space of such families by ${\cal J}$.

Let $\X={\rm Map}(\Sigma,M;A)$ be the space of all smooth maps
$u:\Sigma \to M$ which represent the homology class $A\in H_2(M)$ i.e. such
that $u_*([\Sigma])=A\in H_2(M)$, where $[\Sigma]$ denotes the fundamental
class of the surface $\Sigma$ determined the orientation associated to the
complex structure $j$. For simplicity in this paper we will consider
only homology with integer coefficients  and denote it $H_*(M)$.
So in this situation we will use notation $[u]=A\in H_2(M)$.

We shall denote by ${\X}^{1,p}$ completion of
the space $\X$ with respect  to the Sobolev norm $W^{1,p}$. More precisely,
${\X}^{1,p}$ is the space of maps $u:\Sigma \to M$ whose first covariant
derivatives with respect to Riemannian metric on $M$ are of class $L^p$.
Since the manifold $M$ is compact the topology of this norm does not depend
on the choice of the Riemannian metric. In order for the space ${\X}^{1,p}$
to be well-defined we must assume that $p>2$. The tangent space $T_u\X^{1,p}$
of $\X^{1,p}$ at a smooth $u$ is the completion of the space
$C^{\infty}(u^*TM)$ of all smooth sections $\widehat{u}\in \Gamma (u^*TM)$
in the Sobolev norm.

For a family of almost complex structures $J \in {\cal J}$ let us consider
the infinite dimensional vector bundle
${\cal E} \to \X^{1,p} $  where the fiber at $u$ is the space
${\cal E}_u =L^p(\Sigma,\Omega^{0,1}\otimes_J u^*TM)$ of $L^p$-section
of the vector bundle over $\Sigma$ whose fiber over a point $z \in \Sigma$
is the space of $\C$-linear, with respect to $J(z,u(z))$, maps from
$T^{0,1}_z \Sigma$ to $(u^*TM)_z$.

We note that the zero set of the section
  \[
     {\overline{\partial}}_J: \X^{1,p} \to {\cal E}
  \]
of the infinite dimensional vector bundle  ${\cal E} \to \X^{1,p} $
given by the formula
  \[
     {\overline{\partial}}_J(u)= du + J \circ du \circ j
  \]
is the set of all  $J$-holomorphic curves in the class $A\in H_2(M)$.

We shall denote by $\Omega^{0,1}_{\Sigma}\otimes_J TM$ the vector bundle
over the space $\Sigma \times M$ whose fiber over a point
$(z,m) \in \Sigma \times M$ is the space of $\C$-antilinear maps from
$T_z \Sigma$ to $(TM)_m$ with respect to $J(z,m)$.

Now we want to introduce the general setting for introducing and
dealing with the concept of perturbed $J$-holomorphic curves. Let
$\Omega^1 \otimes \C$ denote the space of all complex valued one
form defined over the product $\Sigma \times M$. On the product
$\Sigma \times M$ there is the almost-complex structure $j \times
J$. With respect to this almost-complex structure we have the
following decomposition of the space $\Sigma \times M$ into the
direct sum
  \[
   \Omega^1 \otimes \C = \Omega^{1,0} \oplus \Omega^{0,1}
  \]
complex linear one forms and complex antilinear one forms.

Any one form in $\Omega^{1,0}$ can be uniquely written in the form
  \[
   \alpha - i \alpha \circ (j \times J),
  \]
where $\alpha$ is a real valued one form in $\Omega^1$.
Similarly, any one form in $\Omega^{0,1}$ can be uniquely written as
  \[
   \alpha + i \alpha \circ (j \times J).
  \]
Our interest will be in following subspaces:
  \begin{enumerate}
   \item
    $\Omega^{1,0}_M$ = subspace of all complex linear one form in
    $\Omega^{1,0}$ whose restriction to the horizontal subbundle
    $T_{\Sigma}(\Sigma \times M)$ of the tangent bundle
    $T(\Sigma \times M)$ is trivial.
   \item
    $\Omega^{0,1}_{\Sigma}$ =  subspace of all complex linear one form in
    $\Omega^{0,1}$ whose restriction to the vertical subbundle
    $T_M(\Sigma \times M)$ of the tangent bundle
    $T(\Sigma \times M)$ is trivial.
  \end{enumerate}
We will use the $\overline{\partial}_{\Sigma}$ operator in the direction
of ${\Sigma}$ in three different situations:
  \begin{enumerate}
   \item
    \[
     \overline{\partial}_{\Sigma}:C^{\infty}(\Sigma \times M)\to
     \Omega^{0,1}_{\Sigma},
    \]
    \[
     \overline{\partial}_{\Sigma}(f)=d_{\Sigma}f+i(d_{\Sigma}f)\circ j.
    \]
   \item
    \[
     \overline{\partial}_{\Sigma}:C^{\infty}(\Sigma \times M,\pi^* TM)
      \to C^{\infty}(\Sigma \times M,\Omega^{0,1}_{\Sigma}\otimes_J TM),
    \]
    \[
     \overline{\partial}_{\Sigma}(X)=d_{\Sigma}X+J\circ (d_{\Sigma}X)\circ j.
    \]
   Here $d_{\Sigma}$ denote the partial derivative in the direction of
   $\Sigma$. Note that it make sense to define $d_{\Sigma}$ since all
   tangent spaces of the form $(\pi^* TM)_{(z,m)}$ are canonically
   identified with $TM_m$.
  \item
    \[
     \overline{\partial}_{\Sigma}:\Omega^{1,0}_M \to \Omega^{0,1}_{\Sigma}
     \otimes_{\C}\Omega^{1,0}_M,
    \]
    \[
     \overline{\partial}_{\Sigma}(\psi-i\psi\circ J)=
     \overline{\partial}_{\Sigma}(\psi)-i\left(
     \overline{\partial}_{\Sigma}(\psi)\right)\circ J,
    \]
   where
    \[
     \overline{\partial}_{\Sigma}(\psi)=d_{\Sigma}(\psi)+i\left(d_{\Sigma}
     (\psi)\right)\circ j
    \]
   for any real valued one form $\psi$ which is zero on $T_{\Sigma}
   (\Sigma \times M)$.
  \end{enumerate}
   We will also need the following maps:
  \begin{enumerate}
   \item
    \[
    \partial_M:C^{\infty}(\Sigma \times M)\to \Omega^{1,0}_M,
    \]
    \[
    \partial_M(f)=d_Mf-i(d_Mf)\circ J.
    \]
   \item
    \[
    \partial_M: \Omega^{0,1}_{\Sigma} \to \Omega^{0,1}_{\Sigma}\otimes
    \Omega^{1,0}_M,
    \]
    \[
    \partial_M(\psi+i\psi \circ j)=\partial_M(\psi)+i\left(\partial_M(\psi)
    \right)\circ j,
    \]
   where
    \[
    \partial_M(\psi)=d_M(\psi)-i\left(d_M(\psi)\right)\circ J
    \]
   for any real valued one form $\psi$ which vanishes on
   $T_M(\Sigma \times M)$.
  \end{enumerate}

  We will often identify the tangent bundle $\pi^*(TM)$ with that of
$\Omega^{1,0}_{M}$ using the metric corresponding to $J$. The
identification is given by the formula
  \[
   X\mapsto \Phi(X),
  \]
where
\begin{equation}\label{identification}
   \Phi(X)(Y)=\omega(X,Y)-i\omega(X,JY)=\inner{X}{JY}_J+i\inner{X}{Y}_J
 \end{equation}
for any $X,Y \in \pi^*(TM)_{(z,m)}$.

Note that the map $\Phi$ is complex antilinear as it should be
since the bundle $\Omega^{1,0}_{M}$ of complex linear forms is a
complex dual of the bundle $\pi^*(TM)$.

Under this identification the following diagram commutes
  \begin{equation}\label{diagram0}
  \def\nornalbaselines{\baselineskip20pt
      \lineskip3pt   \lineskiplimit3pt}
  \def\mapright#1{\smash{
      \mathop{\longrightarrow}\limits^{#1}}}
  \def\mapdown#1{\Big\downarrow\rlap
      {$\vcenter{\hbox{$\scriptstyle#1$}}$}}
   \matrix{C^{\infty}(\pi^*(TM)) &\mapright{\overline{\partial}_{\Sigma}}&
   C^{\infty}(\Omega^{0,1}_{\Sigma}\otimes_J TM) \cr
   \mapdown{\Phi}&          &\mapdown{id \otimes \Phi} \cr
   \Omega^{1,0}_M&\mapright{\overline{\partial}_{\Sigma}}&
   \Omega^{0,1}_{\Sigma}\otimes \Omega^{1,0}_{M}.\cr}
  \end{equation}
  This is because
   \begin{eqnarray*}
    {}\overline{\partial}_{\Sigma}(X)\circ \Phi&=&
      \overline{\partial}_{\Sigma}\omega(X,.)-i\left(
      \overline{\partial}_{\Sigma}\omega(X,.)\right)\circ J\\
   {}&=&d_{\Sigma}\omega(X,.)+id_{\Sigma}\omega(X,.)\circ j\\
   {}&&-i(d_{\Sigma}\omega(X,.))\circ J +d_{\Sigma}\omega(X,.)\circ j \circ J\\
   {}&=&\omega(d_{\Sigma}X,.)-i\omega(d_{\Sigma}X,.)\circ J \\
   {}&& +i(\omega(d_{\Sigma}X,.)-i\omega(d_{\Sigma}X,.)\circ J )\circ j\\
   {}&=&\Phi(d_{\Sigma}X)+i\Phi(d_{\Sigma}X)\circ j\\
   {}&=&id\otimes \Phi(d_{\Sigma}X+J\circ d_{\Sigma}X\circ j)\\
   {}&=&(id\otimes \Phi)\circ \overline{\partial}_{\Sigma}(X).
  \end{eqnarray*}

  \begin{definition}
   An element $ P \in C^{\infty}(\Sigma \times M, \Omega^{0,1}_{\Sigma}
   \otimes_J TM)$ is said to be \textbf{exact} if there is a vector field $X \in
   C^{\infty}(\Sigma \times M, \pi^*TM)$ of the form
    \begin{equation}\label{exact_field}
      X={\Phi}^{-1}\circ \partial_M f
    \end{equation}
   such that
    \begin{equation} \label{exact_term}
     P=\overline{\partial}_{\Sigma}X=d_{\Sigma}+J\circ\left(d_{\Sigma}X
     \right)\circ j.
    \end{equation}
  \end{definition}
  \begin{definition} \label{percurve}
    Let
     $
       \overline{\partial}_{J,f}:\X \to {\cal E}
     $
    be a section of the form
     \begin{equation}\label{exact_element}
       \overline{\partial}_{J,f}(u)= du + J \circ du \circ j + P(u),
     \end{equation}
    where $ P \in C^{\infty}(\Sigma \times M, \Omega^{0,1}_{\Sigma}
    \otimes_J TM)$ is exact, $P=\overline{\partial}_{\Sigma}\circ
    {\Phi}^{-1}\circ \partial_M f$, and
    $P(u)(z)=  P(z, u(z))$.

     The equation
     \begin{equation} \label{PDEH}
        \overline{\partial}_{J,f}(u)=0
     \end{equation}
    will be called a \textbf{ perturbed Cauchy-Riemann}  equation. Solutions
    of the equation (\ref{PDEH}) will be called \textbf{ perturbed
        $J$-holomorphic curves}.
   \end{definition}

  \begin{definition}\label{real_percurve}
    The equation (\ref{PDEH}) is said to be a \textbf{properly perturbed
    Cauchy-Riemann} equation if there is a constant $z \in \mathbb{C}$
    such that the function $f$ satisfies the equation
    \begin{equation}\label{proper_function}
      f=z\,g
    \end{equation}
   for some real function
   \begin{math}
     g \in C^{\infty}(\Sigma \times M, \mathbb{R}).
   \end{math}

      The perturbation term (\ref{exact_term}) is said to be \textbf{properly exact}
     if the function $f$ in the equation (\ref{exact_field})satisfies 
     the equation (\ref{proper_function})
  \end{definition}
\section{Properties of Perturbed $J$-holomorphic Curves}

 \begin{lemma}\label{exact}
  Let $P=\overline{\partial}_{\Sigma}\circ{\Phi}^{-1}\circ \partial_M f$
  be an exact perturbation. Then it can be written as
   \begin{equation}\label{tree}
     P=(id \otimes \Phi)^{-1}\circ \partial_M \circ
       \overline{\partial}_{\Sigma}f.
   \end{equation}
  On the other hand, any perturbation of the form (\ref{tree}) is
  exact.
  \end{lemma}
 \proof{}
  Let us first verify that the following diagram commutes
\begin{equation}
\label{diagram1}
  \def\nornalbaselines{\baselineskip20pt
      \lineskip3pt   \lineskiplimit3pt}
  \def\mapright#1{\smash{
      \mathop{\longrightarrow}\limits^{#1}}}
  \def\mapdown#1{\Big\downarrow\rlap
      {$\vcenter{\hbox{$\scriptstyle#1$}}$}}
   \matrix{C^{\infty}(\Sigma \times M)
   &\mapright{\overline{\partial}_{\Sigma}}&
   \Omega^{0,1}_{\Sigma}\cr
   \mapdown{\partial_M}&          &\mapdown{\partial_M} \cr
   \Omega^{1,0}_M &\mapright{\overline{\partial}_{\Sigma}}&
   \Omega^{0,1}_{\Sigma}\otimes \Omega^{1,0}_{M}.\cr }
 \end{equation}
  Indeed, we may assume that the function $f$ is real. Then
   \begin{eqnarray*}
    {}\partial_M(\overline{\partial}_{\Sigma}f)
    {} &=&\partial_M(d_{\Sigma}f+i(d_{\Sigma}f)\circ j)\\
    {} &=&\partial_M(d_{\Sigma}f)+i(\partial_M(d_{\Sigma}f))\circ j\\
    {} &=&d_M(d_{\Sigma}f)-i(d_M(d_{\Sigma}f))\circ J\\
    {}&& +i(d_M(d_{\Sigma}f))\circ j +(d_M(d_{\Sigma}f))\circ J \circ j
   \end{eqnarray*}
  and
   \begin{eqnarray*}
    {}\overline{\partial}_{\Sigma}(\partial_Mf)
  {}&=&\overline{\partial}_{\Sigma}(\partial_Mf-i(\partial_Mf)\circ J)\\
  {}&=&d_{\Sigma}(d_Mf)+i(d_{\Sigma}(d_M))\circ j\\
  {}&&- i(d_{\Sigma}(d_M))\circ J + (d_{\Sigma}(d_Mf))\circ j \circ J.
   \end{eqnarray*}
  Since partial derivatives commutes the commutativity of the diagram follows.
  To finish the proof of the lemma combine the two commutative diagrams
  (\ref{diagram0}) and (\ref{diagram1}).
  \begin{proposition}\label{exact1}
   Let $\Sigma = S^2$ be the Riemannian sphere and let $\psi \in
   \Omega^{0,1}_{\Sigma}.$ Then the perturbation
    \[
     (id \otimes \Phi)^{-1}\circ \partial_M \psi
    \]
   is exact.
  \end{proposition}
  \proof{}
  Using the Lemma \ref{exact} it is enough to show that $\psi=
 \overline{\partial}_{\Sigma}(f)$ for some function $f \in
  C^{\infty}(\Sigma \times M)$. For every $m \in M$ we have
  $\psi_m = \psi(.,m)\in \Lambda^{0,1}(S^2)$. Since $\Lambda^{0,2}(S^2)=0$,
  $\overline{\partial}\psi_m=0$. Moreover, since $H^{0,1}(S^2)=0$ and
  $H^0(S^2)=\C$ there is the unique function $f_m :S^2 \to \C$ such that
  $\overline{\partial}f_m=\psi_m$ and $f_m(z_0)=0$ for the fixed point
  $z_0 \in S^2$. Define $f$ by the formula
   \[
     f(z,m)=f_m(z).
   \]
  Then one verifies that $\psi= \overline{\partial}_{\Sigma}(f)$ and this
  finishes the proof of the Proposition. 

 We note that if $g: \mathbb{R} \times S^1 \to S^2$, $g(s,t)=z=s+it$ is a
 holomorphic coordinate
 system and a function $H: \mathbb{R} \times S^1 \times M \to \mathbb{R}$ has 
 a compact support then the Proposition \ref{exact1} implies that the perturbation
  \begin{equation}\label{hamiltonian_perturbation}
    P=\nabla H ds - J \nabla H dt
  \end{equation}
 is exact. This is so because
  \begin{eqnarray*}
   {}&&id \otimes \Phi (\nabla H ds -J \nabla H dt)\\
   {}&=&(idH+dH \circ J) \otimes (ds-idt)\\
   {}&=&i(dH-idH \circ J) \otimes ds +(dH-idH \circ J)\otimes dt\\
   {}&=&\partial_M (iHds +H dt).
  \end{eqnarray*}
 In our construction of a new Floer symplectic (co)homology it will be very important
 to note that the Hamiltonian perturbation term (\ref{hamiltonian_perturbation})
 is in fact properly exact. It follows from the equation
  \begin{equation}\label{hamiltonian_exactly}
    \overline{\partial}_{\Sigma}(ig)=iH\,ds+H\,dt,
  \end{equation}
 where the function $g$ is defined as
  \[
  g(z,m)=\int_{(\infty,m)}^{(z,m)}\,H\,ds + H\,dt,
  \]
 where the integration is taken over any smooth path contained in the set
 $\Sigma \times m$
 \begin{theorem}\label{noholomorphic}n
  Let $u:\Sigma \to M$ be a smooth map from a Riemannian surface
  $\Sigma$ to a symplectic manifold $M$ with compatible family of
  almost-complex structures $J$. Assume that $du(z_0)\not=0$ and
  that $du(z_0)$ is complex anti-linear at a single point $z_0 \in
  \Sigma$.

  Then $u$ is not a perturbed $J$-holomorphic curve.
 \end{theorem}
 \proof{}
  Assume that $u$ is a perturbed $J$-holomorphic curve. Then by
  the Definition \ref{percurve} and the Lemma \ref{exact} it
  satisfies the equation
    \[
      du + J \circ du \circ j + P(u)=0,
    \]
  where
\begin{equation}\label{fcurve}
  P=(id \otimes \Phi)^{-1}\circ \partial_M \circ
       \overline{\partial}_{\Sigma}f.
\end{equation}
  Consider the tensor $(id \otimes \Phi)\circ P \in
  \Omega^{0,1}_{\Sigma}\otimes \Omega^{1,0}_{M}$.
  For any two vectors $v,w \in T_{z_0} \Sigma$ we have the
  following formula
  \begin{eqnarray*}
    {}\theta(v,u)&:=&(\textrm{id} \otimes \Phi)\circ P(v\otimes
    du(w))\\
    {}&=&-\frac{1}{2} \inner{\overline{\partial}_Ju(v)}{J\circ
    \overline{\partial}_Ju(w)} -\frac{1}{2}i
    \inner{\overline{\partial}_Ju(v)}{\overline{\partial}_Ju(w)}.
   \end{eqnarray*}
  This is because of the formula (\ref{identification}) and the
  fact that, since $du(z_0)$ is complex anti-linear, at $z_0$ we
  have
   \[
     \overline{\partial}_Ju(z_0)=-P(u)(z_0)=2du(z_0).
   \]
  In particular, since vectors $\overline{\partial}_Ju(v)$ and
  $J\circ \overline{\partial}_Ju(v)$ are orthogonal we obtain that
  the quadratic function
   \[
    \theta(v,v)=-\frac{1}{2}i
    \inner{\overline{\partial}_Ju(v)}{\overline{\partial}_Ju(v)}
   \]
  is purely imaginary. Moreover, one easily verifies that
\begin{equation}\label{jinvariant}
  \theta(j\circ v,j\circ v)=\theta(v,v).
\end{equation}
  Therefore, because of (\ref{fcurve}) the expression $\partial_M \circ
       \overline{\partial}_{\Sigma}f(v,du(v))$ is also purely
  imaginary at $z_0$. Thus one computes at $z_0$
  \begin{eqnarray*}
  {}\partial_M \circ
       \overline{\partial}_{\Sigma}f(v,du(v))
    &=&i(d_{jv}d_{du(v)}f_1+d_vd_{du(jv)}f_1)\\
  {}&&+i(d_vd_{du(v)}f_2-d_{jv}d_{du(jv)}f_2,
  \end{eqnarray*}
  where we have written $f=f_1+if_2$.

  Now, one can easily see that at $z_0$
   \[
     \partial_M \circ \overline{\partial}_{\Sigma}f(v,du(v))=-
     \partial_M \circ \overline{\partial}_{\Sigma}f(j\circ v,
     du(j\circ v)).
   \]
  This is a contradiction because of the equation
  (\ref{jinvariant}). This proves that $u$ can not be a perturbed
  $J$-holomorphic curve. 
\begin{remark}{}
 The Theorem \ref{noholomorphic} implies that the antipodal map
 $u:S^2 \to S^2$ where $S^2$ is equipped with the standard complex
 structure is not a perturbed $J$-holomorphic curve.
\end{remark}
\section{Hermitian Structures}
Here we will review basic properties of Hermitian structures on
almost-complex manifolds and apply them in our context. Apart from
other sections of these notes we will consider almost-complex
manifold $(M,J)$ with a fixed almost-complex structure $J$.

Let $E \to M$ be a complex vector bundle over $M$. Then the
almost-complex structure $J$ induces the following decomposition
  \[
   \Omega^k(M,E)=\bigoplus_{p+q=k}\Omega^{p,q}(M,E),
  \]
where $\Omega^k(M,E)$ denotes the space $\Omega^k(M)\otimes_{\C}E$
of smooth $E$-valued k-forms and $\Omega^{p,q}(M,E)$ denotes its
subset of all k-form complex linear with respect to p arguments
and complex anti-linear with respect q arguments.

Let $\nabla:C^{\infty}(M,E)\to \Omega^1(M,E)$ be a covariant
derivative on the vector bundle $E$. We can decompose $\nabla$ as
\begin{equation}\label{decompose}
 \nabla=\partial_{\nabla}+\overline{\partial}_{\nabla},
\end{equation}
into complex linear part and complex anti-linear part,
respectively.

The complex linear part $\partial_{\nabla}:C^{\infty}(M,E)\to
\Omega^{1,0}(M,E)$ is given by the formula
\begin{equation}\label{complex linear}
  \partial_{\nabla}=\frac{1}{2}\left(\nabla-i\,\nabla\circ
  J \right)
\end{equation}
and the complex anti-linear part
$\overline{\partial}_{\nabla}:C^{\infty}(M,E)\to
\Omega^{0,1}(M,E)$ is given by the formula
\begin{equation}\label{complex anti-linear}
  \overline{\partial}_{\nabla}=\frac{1}{2}\left(\nabla+
  i\,\nabla\circ
  J \right).
\end{equation}

\begin{definition}\label{Cauchy-Riemann operator}
 An operator $D'':C^{\infty}(M,E)\to\Omega^{0,1}(M,E)$ is said to
 be a Cauchy-Riemann operator if
\begin{equation}\label{Cauchy-Riemann equation}
  D''(fs)=\overline{\partial}\,f \otimes s + f\,D''s,
\end{equation}
 for every $f\in C^{\infty}(M,\C)$ and any $s \in
 C^{\infty}(M,E)$.
\end{definition}

\begin{lemma}\label{trivial}
 The complex anti-linear part $\overline{\partial}_{\nabla}$ of
 the covariant derivative $\nabla$ is a Cauchy-Riemann
 operator on $E$.
\end{lemma}
\begin{definition}\label{Hermitian metric}
  A Hermitian metric $h$ in $E$ is a smooth family of Hermitian
  inner products in the fibers of the vector bundle $E$.
\end{definition}

As a example, let $(M,J)$ be a symplectic manifold and let $J$ be
$\omega$-compatible almost-complex structure on $M$. Then the
tensor $\inner{\cdot}{\cdot}$ given by the formula
\begin{equation}\label{inner Hermitian}
  \inner{u}{v}=\omega(Ju,v)+i\,\omega(u,v)
\end{equation}
defines a Hermitian metric on the tangent bundle $TM$ with respect
to the almost-complex structure $J$. Compare this with the formula
(\ref{identification}).

Any Hermitian metric $h$ on $E$ determines the Hermitian
connection $\nabla_h=\nabla$ on $E$. It is a unique connection
  \[
    \nabla :C^{\infty}(M,E)\to\Omega^1(M,E)
  \]
which preserves both the complex structure of $E$ and the metric
$Re\,h$ ( the real part of $h$) induced by the Hermitian structure
$h$.

Alternatively, the Hermitian connection is a unique connection
$\nabla$ such that
\begin{equation}\label{alternative connection}
  d\,h(u,v)=h(\nabla\,u,v)+h(u,\nabla\,v).
\end{equation}
Here is the basic fact about Hermitian connections:

\begin{proposition}\label{every}
  For every Cauchy-Riemann operator $D'':C^{\infty}(M,E)\to
  \Omega^{0,1}(M,E)$ there exists a unique Hermitian connection
  $\nabla$ such that its complex anti-linear part
    \[
      \overline{\partial}_h=\frac{1}{2}(\nabla+i\,\nabla \circ J)
    \]
  is equal to $D''$.
\end{proposition} 

 A Hermitian metric $h$ on a complex bundle $E$ induces the
 Hermitian metric $h^*$ on the complex dual bundle $E^*$. Namely,
 If $e=(e_i)$ is a unitary frame for $E$, $e^*=(e^*_i)$ the dual
 frame for $E^*$, then set
   \[
     h^*(e^*_i,e^*_j)=\delta_{ij}.
   \]
 If we identify $E$ with $E^*$ ( in complex anti-linear way) via
 formula
\begin{equation}\label{via formula}
  s \mapsto s^*=h(\cdot,s)
\end{equation}
 then the formula (\ref{alternative connection}) implies that the
 Hermitian connection $\nabla^*$ on the dual bundle $E^*$ is
 uniquely determined by the requirement:
\begin{equation}\label{requirement}
  d\,\inner{t}{s}=\inner{\nabla^*\,t}{s}+\inner{t}{\nabla\,s }
\end{equation}
 for $t\in C^{\infty}(M,E^*)$ and $s\in C^{\infty}(M,E)$.

 This is so because under the identification (\ref{via formula})
 the connection $\nabla$ corresponds to a connection which is both
 preserving the metric and the complex structure on $E^*$ and thus
 it corresponds to $\nabla^*$
  \begin{theorem}\label{for and}
   Let $(M,J)$ be an almost-complex manifold and let $h$ be a
   Hermitian metric on the tangent bundle $E=TM$. Then the complex
   anti-linear part
   \[
    \overline{\partial}_{h^*}:\Omega^{1,0}(M) \to
    \Omega^{0,1}\left(\Omega^{1,0}(M)\right)
    \cong \Omega^{0,1}(M)\otimes \Omega^{1,0}(M)
   \]
   of the dual Hermitian connection $\nabla^*$ is given by the following 
  formula
   \begin{equation}\label{antiholomorphic_part}
    \iota_{X}(\overline{\partial}\,\eta)=
      -\frac{1}{2}\left(\iota_{X}(d\,\eta) +i\,\iota_{J\,X}(d\,\eta)\right)
   \end{equation}
  for any vector field $X$.
   \end{theorem}
   Note that we have identified the tangent space $TM$ with bundle
   $\Omega^{1,0}(M)$ of complex linear forms via formula
   (\ref{via formula})
  \proof{} The formula (\ref{requirement}) implies that for any
  complex linear form $\eta \in \Omega^{1,0}(M)$ and for any
  vector field $X$ we have
\begin{equation}\label{we have}
  d(\iota_X\,\eta)=\iota_X(\nabla^*\eta)+\eta(\nabla X).
\end{equation}
  Let us fixed an arbitrary point $m \in M$ and let $X_m \in T_m(M)$
  and $Y_m \in T_m(M)$ be any two tangent vectors at $m$. Choose
  a smooth map
   \[
    \sigma : (-\epsilon,\epsilon)^3 \to M
   \]
  with $\sigma(0,0,0)=m$ satisfying the following three
  conditions:
\begin{enumerate}
  \item $d\,\sigma_{(0,0,0)}([1,0,0])=X_m$,\,\,\,\,
        $d\,\sigma_{(0,0,0)}([0,1,0])=J\,X_m$,
  \item $d\,\sigma_{(0,0,0)}([0,0,1])=Y_m$,
  \item $\eta(d\,\sigma_{(s,t,0)}([0,0,1]))= const$.
\end{enumerate}
  Let $X_1$, $X_2$, and $Y$ denote vector fields
  $d\,\sigma ([1,0,0])$, $d\,\sigma ([0,1,0])$, and \\
  $d\,\sigma ([0,0,1])$, respectively. With this notation  the condition 3. above implies
  that
   \[
     d(\iota_Y\,\eta)(X_m)=d(\iota_Y\,\eta)(J\,X_m)=0.
   \]
  Therefore combining general identities for exterior derivative
   \[
    \iota_{X_i}(d\,\eta)(Y)=d\,\eta(X_i,Y)=d(\iota_Y\,\eta)(X_i)-
    d(\iota_{X_i}\,\eta )(Y)-\eta ([X_i,Y])
   \]
  with the identity (\ref{we have})  and noting that
  commutators $[X_i,Y]$ are trivial we obtain
   \[
    \iota_{X_i}(d\,\eta)(Y)=-\iota_{X_i}(\nabla^*\,\eta)(Y)-
    \eta(\nabla\,X_i)(Y)
   \]
  for any tangent vector $Y$ at $m$.

  Thus at the point $m$ we have
\begin{equation}\label{first}
     \iota_{X_m}(d\,\eta)=-\iota_{X_m}(\nabla^*\,\eta)-\eta(\nabla\,X_m)
\end{equation}
  and
\begin{equation}\label{second}
   \iota_{J\,X_m}(d\,\eta)=-\iota_{J\,X_m}(\nabla^*\,\eta)-
   \eta(\nabla\,(J\,X_m)).
\end{equation}
  Now notice that the expression  $\eta(\nabla\,X)$ is complex linear
  with respect to the variable $X$ since the metric $\nabla$ is
  Hermitian and the 1-form $\eta$ is chosen to be complex linear.
  Therefore  combining identities
  (\ref{first}) and (\ref{second}),  
  at the point $m \in M$ we obtain the equation
   \[
     \iota_{X_m}(\overline{\partial}\,\eta)=
      -\frac{1}{2}\left(i_{X_m}(d\,\eta) +i\,\iota_{J\,X_m}(d\,\eta)\right).
   \]
  Since this equality is true for any point $m\in M$ and any
  tangent vector $X_m \in T_m M$ the theorem follows.

 %
 %
\begin{proposition}\label{gradient}
Let $(M,\omega)$ be a symplectic manifold of the dimension $2n$
with a compatible almost-complex structure $J$. Denote by $h'$  a
Hermitian metric on the tangle bundle $TM$ by the equation
(\ref{inner Hermitian}) and by $\,\Phi:TM \cong \Omega^{1,0}(M)$
the isomorphism given by the equation (\ref{identification}). Let
$f:M \to \mathbb{C}$ be a complex function of the form $f=z \zeta$ for some
complex number $z$ and a real valued function $\zeta$ such that
$\partial\,(f)(m) \not= 0$ for $m \in M$. Then there is an open
set $U$ with $m \in U$ and a Hermitian metric $h$ on the set $U$
with the properties:
\begin{itemize}
 \item
  $h(X_1,X_2)=h'(X_1,X_2)$ for any two vector fields $X_i$,
  $i=1,2$ with $h'(X_i, \Phi^{-1} \circ \partial\,(f))\cong 0$ on $U$. Moreover,
  $h'(X, \Phi^{-1} \circ \partial\,(f))=0$ if and only if $h(X, \Phi^{-1} 
  \circ \partial\,(f))=0$ for any vector field $X$ on $U$.
 \item
   The map
     \begin{eqnarray*}
       {}U & \to & \mathbb{R}\\
       {}u &\mapsto& h(\Phi^{-1} \circ \partial \,(f), \Phi^{-1} \circ
                  \partial\,(f))\\
       {}   && =h(\Phi^{-1} \circ \partial \,(f)(m), \Phi^{-1} \circ
                  \partial\,(f)(m)) 
     \end{eqnarray*}
    is constant.
 \item
   $\nabla(\Phi^{-1} \circ \partial\,(f)) \in TM \otimes_J \Omega^{1,0}(U)$,
   where $\nabla$ denotes the Hermitian connection of the metric
   $h$.
\end{itemize}
\end{proposition}
\proof{}
We first note that the first two conditions determine the metric $h$ uniquely.
Thus we only need to show that $\nabla(\Phi^{-1} \circ \partial\,(f)) $ 
is complex linear. Without loss of generality we may assume that $f$ is 
a purely imaginary function: $f=i\zeta$. Next,
we note that
 \[
   h'(\,\cdot\,,\Phi^{-1}\circ \partial(f))= \partial \zeta.
 \]
Let 
\[
  \eta:=h(\,\cdot\,,\Phi^{-1}\circ \partial(f)).
\]
By the construction the forms $\partial \zeta$ and $\eta$ differ only by a factor
of real function and thus have the same zero sets.
Using the Theorem \ref{for and} it is enough to show that
\[
 d \,(\eta)=0.
\]
Define function $g:M\to \mathbb{R}$ 
 \[
 g(x)=
   \left\{
       \begin{array}{ll}
         0 & \mbox{if $x \in \zeta^{-1}(m)$}\\
         t & \mbox{where $t$ is a time needed to travel}\\ 
         {}&\mbox{from the level $\zeta^{-1}(m)$ to $x$ along the flow of 
            $\Phi^{-1} \circ \partial\,(f)$.}      
       \end{array}\right.
\]
By construction of the form $\eta$ if $\eta(Y)=0$ then the vector $Y$ is tangent
to a level set of the function $g$. Now since $g$ is a real function we have
\[
  \oint_{\gamma}\partial g=0
\] 
over any closed (contractible) loop $\gamma$. This implies that 
the function $\rho :M \to \mathbb{C}$ 
\[
   \rho(x)=\int_m^x \partial g
\]
is well defined in a (contractible) neighborhood of $m$. 
Moreover simple computation shows that
\[
  d \rho = \eta.
\]
This finishes the proof of the Proposition since we have now
\[
  d \eta =d\,d \rho=0.
\]

 For vector field whose connection is of type
 $(1,0)$ we have the following:
\begin{lemma}\label{torsion is trivial}
 Let $(M,J)$ be an almost-complex manifold of dimension $2n$ and
 let $h$ be a Hermitian metric defined on the tangle bundle $TM$.
 Let $X$ be a vector field on $M$. Assume that
 \[
  \nabla(X) \in TM \otimes \Omega^{1,0}(M).
 \]
 Then for any point $m \in M$ and any tangent vector $Y_m \in T_m(M)$
 the torsion $T(X,Y_m)$ is trivial.
\end{lemma}
\proof{}
 Choose a vector field $Y$ which agree with the tangent vector
 $Y_p$ at the point $m$ such that commutators $[X,Y]$ and $[JX,Y]$
 are trivial. Using the fact that the Torsion tensor $T(X,Y)$ is
 complex anti-linear with respect to the two variables $X$ and
 $Y$ (see the Section \ref{linearization section}) we have
 \begin{eqnarray*}
  {}JT(X,Y_m)&=&-T(JX,Y)\\
        {}&=&-\nabla_{JX}(Y) + \nabla_{Y}(JX)\\
        {}&=&-J(\nabla_X(Y)-\nabla_Y(X))\\
        {}&=&-JT(X,Y_m).
 \end{eqnarray*}
 This shows that $T(X,Y_m)=0$. 

\section{ Compactness }
Let us consider the following weak version of the Gromov's Compactness
Theorem \cite{G}
    \begin{theorem}
     \label{Gromov}
     Let $(M,\omega)$ be a compact symplectic manifold and let $J_k$
     be a sequence of $\omega$-tame almost complex structures which
     converge to $J_{\infty}$ in $C^{\infty}$-topology. Then for any sequence
     $u_k: \Sigma \to M$ of $J_k$-holomorphic curves with uniformly bounded
     energy there are subsequence (still denoted by $u_k$),
     a finite collection $(u^1,...,u^m)$ of $J_{\infty}$-holomorphic
     spheres $u^i:S^2 \to M$, and a $J_{\infty}$-holomorphic curve
     $u^{\infty}:\Sigma \to M$
     such that for corresponding homology classes in $H_2 (M)$ we have
       \[
           [u_k]=[u^{\infty}]+[u^1]+...+[u^m]
       \]
    for all $k$ large enough.
   \end{theorem}
We want to prove a similar theorem for perturbed $J$-holomorphic curves
satisfying the equation (\ref{PDEH}).  Let us describe an extension of
Gromov's nice trick to perturbed $J$-holomorphic curves.
Consider a solution $u:\Sigma \to M$ of the differential equation
(\ref{PDEH}).  To such $u$ we associate a map
$ \widetilde{u}:\Sigma \to \Sigma \times M$ given by the formula
    \[
          \widetilde{u}(z) = (z, u(z)),
    \]
for $z= (s,t)\in \Sigma$. Then $\widetilde{u}$ satisfies the nonlinear
Cauchy-Riemann equation
    \begin{equation}
     \label{hol}
       \overline{\partial}_{\widetilde{J}}(\widetilde{u})= d\widetilde{u}
       +\widetilde{J} \circ d\widetilde{u} \circ j=0 ,
   \end{equation}
where $\widetilde{J}$ is almost-complex structure on the product
$\Sigma \times M$ given by the formula
    \[
     \widetilde{J}=
       \left[
         \matrix{j                &       0 \cr
                        - Pj &  J \cr}
       \right].
     \]
To check that $\widetilde{u}$ satisfies the equation (\ref{hol}) note that
    \[
     \overline{\partial}_{\widetilde{J}}(\widetilde{u})=
       \left[
         \matrix{ 0   \cr
                 \overline{\partial}_{J}(u)+P(u)     \cr}
       \right].
    \]
Choose a symplectic form $\omega_0$ on $\Sigma$ such that the
complex structure $j$ on $\Sigma$ is compatible with $\omega_0$
and such that $\omega_0([\Sigma])=1$.  Define a symplectic structure
$\widetilde{\omega}$ on $\Sigma \times M$ by the formula
    \[
         \widetilde{\omega}= N\omega_0 +\omega,
     \]
where $N$ is a positive number.
   \begin{lemma}
      If $N$ is large enough, then the almost-complex structure
      $ \widetilde{J}$ is  $ \widetilde{\omega}$-tame.
   \end{lemma}
    \proof{}
  Define
     \[
      f = ||\omega||_{L^{\infty}} \sup_{(z,m) \in \Sigma \times M} ||P(z,m)||
     \]
  and compute
     \begin{eqnarray*}
       {}  \widetilde{\omega} \left( \widetilde{J}
                          \left[
                           \matrix{ a   \cr
                                             v   \cr}
                         \right] ,
                         \left[
                           \matrix{ a   \cr
                                             v   \cr}
                         \right]
                         \right)
                           &=& N\omega_0  ( j(a), a)+
                          \omega(P(a) +Jv, v)\\
        {}  & \geq & N|a|^2-|a|f|v|_{\widetilde{J}}+
                     |v|_{\widetilde{J}}^2\\
        {}  & =      &((N)^{\frac{1}{2}}|a|-|v|)^2+
                                 (2(N)^{\frac{1}{2}}-f)|a||v| \\
        {}  & \geq &(2(N)^{\frac{1}{2}}-f)|a||v| .
   \end{eqnarray*}
  Here we have used the notation
     \[
       |v|^2=|v|_{\widetilde{J}}^2= \sup_{z} \omega(J(z)v,v).
     \]
  To finish the  proof of the lemma it is enough to choose
  $ N>\left(\frac{f}{2}\right)^2$. 

So we choose $N$ large enough so that the almost-complex structure
$ \widetilde{J}$ is  $ \widetilde{\omega}$-tame. In this situation there is
a Riemannian metric defined on the product $\Sigma\times M$ via the formula
   \[
   { \inner{v}{w}}_{\widetilde{J}}=\frac{1}{2} \left(  \widetilde{\omega}
     (\widetilde{J}v, w)+\widetilde{\omega}(\widetilde{J}w,v) \right),
   \]
for tangent vectors $v,w$.

The energy of  $\widetilde{u}$ is defined as
   \[
     E(\widetilde{u})= \frac{1}{2} \int \int \left|D \widetilde{u}
     \right|_{\widetilde{J}}^2.
    \]
In fact the energy depends only on a homology class $[\widetilde{u}]$ of
$\widetilde{u}$ and we have
  \[
    \widetilde{\omega}(\widetilde{u})=E(\widetilde{u})\geq 0.
  \]

    \begin{theorem}
     \label{Gromov2}
     Let $(M,\omega)$ be a compact symplectic manifold and let $J_{\nu}$
     be a sequence in ${\cal J}$ of $\omega$-compatible families of
      almost complex  structures which converge to $J_{\infty}$
      in $C^1$-topology.
     Let $ P_{\nu}$ be a sequence
     of perturbations of the form
       \[
         P_{\nu} \in C^{\infty}(\Sigma \times M,
         \Omega^{0,1}_{\Sigma}\otimes_{J_{\nu}} TM)
       \]
      which converge to
     $ P_{\infty}$ in $C^{1}$-topology. Then for any sequence
     $u_{\nu}: \Sigma  \to M$ of solutions of the equation (\ref{PDEH})  with
      $P=P_{\nu}$, $J=J_{\nu}$ such that  $[u_{\nu}]=A\in H_2(M)$
      there is subsequence (still denoted by $u_{\nu}$) and a finite
     collection
     $(u_{\infty};C_1...,C_m)$, where $u_{\infty}$ is a of solution of the
       equation (\ref{PDEH})  with
      $P=P_{\infty}$ and $J=J_{\infty}$ and $C_i$ are
     $J_{\infty}(z_i)$-holomorphic spheres $C_i:S^2 \to M$
     such that for corresponding homology classes in $H_2 (M)$ we have
      \begin{equation}
         \label{A}
           A=[u_{\infty}]+[C_1]+...+[C_m].
       \end{equation}
    \end{theorem}
  \proof{}
  Consider corresponding  elements ${\widetilde{u}}_{\nu}$,
  ${\widetilde{J}}_{\nu}$, and ${\widetilde{J}}_{\infty}$.
  Since $J_{\nu} \to J_{\infty}$, $P_{\nu} \to P_{\infty}$, and
  all $u_{\nu}$ represent the same homology class $A$ the energy
  of the sequence ${\widetilde{u}}_{\nu}$
  is uniformly bounded for almost all $\nu$ and we may apply
   the Gromov's theorem \ref{Gromov}.
  Thus there is a collection $({\widetilde{u}}^i:S^2 \to \Sigma \times M)$,
  $i=1,...,m$, of ${\widetilde{J}}_{\infty}$-holomorphic spheres
  ${\widetilde{J}}_{\infty}$-holomorphic curve
     ${\widetilde{u}}^{\infty}:\Sigma \to \Sigma \times M$  such that
       \[
       A+[\Sigma]=   [{\widetilde{u}}_{\nu}]=[{\widetilde{u}}^{\infty}]+
        \sum_{i=1}^m [{\widetilde{u}}^i]
       \]
   in $H_2(\Sigma) \oplus H_2 (M) \subseteq H_2 ( \Sigma \times M)$.
    For every $i$, we can write a unique decomposition
      \[
          [{\widetilde{u}}^i]=A^i+B^i,
       \]
    and also a unique decomposition
      \[
          [{\widetilde{u}}^{\infty}]=A^{\infty}+B^{\infty},
       \]
    where $ A^i,A^{\infty} \in H_2 (M)$ and $B^i,B^{\infty} \in H_2 (S^2)$.
    In particular, we have
    \[
    A^{\infty}+ \sum_{i=0}^m A^i = A, \>\>\>\>\>
    B^{\infty}+\sum_{i=0}^m B^i =[\Sigma].
     \]
    Now, the class $[\Sigma]$ is indecomposable in $H_2 (\Sigma)$. Therefore
     $ B^{\infty}=[\Sigma]$ and $B^i=0$.
   Thus the map $ pr_1 \circ {\widetilde{u}}^{\infty} :\Sigma \to \Sigma$,
   where $pr_1: \Sigma \times M \to \Sigma$ is the natural projection on
    the first factor,  is a holomorphic of degree one.
    Eventually reparametrizing the map ${\widetilde{u}}^{\infty}$ we may
    assume  that the map $ pr_1 \circ {\widetilde{u}}^{\infty} $
    is the identity map on $\Sigma$.
     With this reparametrization  define
      \[
           u_{\infty}= pr_2 \circ {\widetilde{u}}^{\infty},
      \]
   where $pr_2: \Sigma \times M \to M$ is the natural projection on the
   second factor. It is easy to verify now, that  $u_{\infty}$ such defined
   satisfies the  equation (\ref{PDEH})  with
      $P=P_{\infty}$ and $J=J_{\infty}$ .

   Consider now maps ${\widetilde{u}}^i$, for $i=1,...,m$.
   The projections  $ pr_1 \circ {\widetilde{u}}^i $ are holomorphic
   maps from $S^2$ to $\Sigma$ of degree zero. Therefore they must
   be constant:
       \[
            pr_1 \circ {\widetilde{u}}^i (S^2)= z_i \in \Sigma.
       \]
    This easily implies that $ pr_2 \circ {\widetilde{u}}^i : S^2 \to M$ is
    $J(z_i)$-holomorphic sphere. Define
       \[
          C_i=pr_2 \circ {\widetilde{u}}^i : S^2 \to M.
       \]
   This finishes the proof of the theorem since the equation (\ref{A})
   is obvious now.
  \begin{remark}
   If $\Sigma$ is a Riemann sphere $S^2$ then the above theorem
   is true if we replace homology classes in $H_2(M)$ by homotopy
   classes in $\pi_2 (M)$.
  \end{remark}
  \begin{remark}
   \label{tilde}
    Note that we have the following identity:
     \[
        \omega(C_i)= \widetilde{\omega}({\widetilde{u}}^i).
     \]
  \end{remark}
\section{Linearization}\label{linearization section}
We will examine the moduli space
  \[
    \M(A, J, P)
  \]
of all solutions of the equation (\ref{PDEH}) which represent a
given homology class  $A\in H^2(M)$. In order to do it we need
study the linearization of the perturbed Cauchy-Riemann operator.

In general, there is no unique way to construct such linearization
since we have to define a way to identify, for each fixed $z \in \Sigma$,
the fiber
  \[
   L^p\left(\Omega^{0,1}(z){\otimes}_{J(z)}T_{u(z)}M\right)
  \]
with that of the form
  \[
    L^p\left(\Omega^{0,1}(z){\otimes}_{J(z)}T_{m}M\right)
  \]
for any $m$ close to $u(z)$.

Perhaps, it will be the best for us if we choose such identification
based on the family of Hermitian connections $\nabla(z)$ relative to the
family of almost complex structures $J(z)$. Recall that a Hermitian
connection is a connection that preserves a Hermitian metric with respect to
the $ J(z)$. In contrast to the Levi-Civita connection its torsion tensor
$T$ is, in general, nontrivial and, moreover, it is complex anti-linear in
two variables, i.e.
   \[
     T(J\xi,\eta)=T(\xi,J\eta)=-JT(\xi,\eta).
   \]

To describe the linearization based on Hermitian connections let us
choose $u \in \M(A, J, P)$
 and $z \in \Sigma$. Since the hermitian connection
$\nabla (z)$ preserves the almost-complex structure $J(z)$ the map
  \[
    \Phi_{\widehat{u}}(z):L^p ((\exp_{u}(\widehat{u}))^*(\Omega^{0,1}
     \otimes_J TM)) \to L^p(\Omega^{0,1} \otimes_J u^*TM)
  \]
induced by the parallel transport along the geodesic curve $t \to
\exp(t\widehat{u})$ corresponding to the Hermitian metric at $z$
on $M$ is well-defined. Thus in the neighborhood of $u$ the
perturbed Cauchy-Riemann $\overline{\partial}_{J,P}$ is
represented by the map
  \[
    {\cal F}:W^{1,p}(u^*TM) \to L^p(\Omega^{0,1} \otimes_J u^*TM)
  \]
defined by
  \[
    {\cal F}(\widehat{u})=\Phi_{\widehat{u}}(z)(\overline{\partial}_{J,P}
    (\exp(\widehat{u}))).
  \]
The linearization at $u$
  \[
    D_u:W^{1,p}(u^*TM) \to L^p(\Omega^{0,1} \otimes_J u^*TM)
  \]
is defined as $D_u(\widehat{u}(z))=d{\cal F}(0)(\widehat{u}(z)).$ It is not
hard to compute $ D_u$ (See \cite{MS} for $J$-holomorphic curves).

 \begin{proposition}
  If $u:\Sigma \to M$ satisfies the equation (\ref{PDEH}) with the exact
  perturbation term $P=\overline{\partial}_{\Sigma}X$ and
  $\widehat{u}\in C^{\infty}(u^*TM)$ then  the operator $D_u(\widehat{u})$ can
  be written as
    \begin{eqnarray*}
    {}  D_u \widehat{u}&=&
    \nabla^*(\widehat{u})+T(du,\widehat{u}) +
     \nabla_{\widehat{u}}(d_{\Sigma}X)\\
     {}&& +J\circ (\nabla^*(\widehat{u})+T(du,\widehat{u})
          + \nabla_{\widehat{u}}(d_{\Sigma}X))\circ j,
     \end{eqnarray*}
  where $\nabla^*$ denotes the induced connection on $u^*(TM)$ via
  the map $(id,u)$.
 \label{D}
 \end{proposition}

  \begin{lemma}
   For $u\in\M(A, J, P)$ the operator $D_u:W^{1,2}(u)\to L^2(u)$
   is elliptic.
   Here $W^{1,2}(u)$ and $L^2(u)$ denote spaces $W^{1,2}(u^*TM)$
   and $W^{0,2}(\Omega^{0,1} \otimes_J u^*TM)$, respectively.
   Thus the operator $D_u$ is Fredholm.
  \label{index}
  \end{lemma}
 \proof{}
 Since the vector bundles  $W^{1,2}(u^*TM)$ and
 $W^{0,2}(\Omega^{0,1} \otimes_J u^*TM)$ are defined
 over the compact manifold $\Sigma$ it is enough to notice that main
 symbol of the operator $F(u)$ is elliptic. From general theory of elliptic
 operators defined on vectors bundles over a compact manifold follows
 that $F(u)$ is also Fredholm.
  \begin{remark}
  \label{fredholm}
    By the elliptic regularity it follows that the map
      \[
        D_u:W^{1,p}(u^*TM) \to L^p(\Omega^{0,1} \otimes_J u^*TM)
      \]
    is also Fredholm.
  \end{remark}
Here is the basic fact about properly exact perturbations.
  \begin{theorem}
   Let the perturbation term $P$ be properly exact, 
   $P=\overline{\partial}_{\Sigma}X$.
    Then there exists a family of Hermitian connections such that
    \[
     D_u((id,u)^*X)
    \]
 is closed to $(id,u)^*(\overline{\partial}_{\Sigma}X)$.
   so $(id,u)^*(\overline{\partial}_{\Sigma}X)$ is in the range of $D_u$.
  \label{exact2}
  \end{theorem}
 \proof{}
 If $X(z,u(z))=0$ then choose  a Hermitian connection $\nabla(z)$ 
 corresponding to the Hermitian metric determined by $(J(z), \omega)$. 
If $||X(z,u(z))||>\epsilon$ then choose the perturbed Hermitian connections
(in the neighborhood of $u(z)$) given by the Proposition \ref{gradient}.
For the latter connections the Proposition \ref{gradient} and the Lemma
\ref{torsion is trivial} imply that 
\[
J \nabla_Y(X)=\nabla_{JY}(X).
\]
Therefore for these connections we have
\[
D_u((id,u)^*X)=(id,u)^*(\overline{\partial}_{\Sigma}X)+
                  \nabla_{du+J\circ du \circ j+P(u)}(X).
\]
The last expression is zero. Now make $\epsilon$ small enough and choose
Hermitian connections for $z$ satisfying $0<||X(z,u(z))||< \epsilon$ so that
we get a smooth bounded (by a constant independent of the choice of
$\epsilon$) family of Hermitian connections satisfying the conclusion of the 
Theorem.

\section{Compactness Properties of the space $V$}

We will study the space of smooth functions (real valued, for simplicity, but 
all applies to the space of function of the form $zf$ where $z$ is a fixed 
complex number and $f$ is arbitrary real valued function) defined over 
the set $\Sigma \times M$.
Choose a decreasing sequence $\epsilon_k>0$ and consider the
subspace $C^{\infty}_{\epsilon}(\Sigma \times M)$ of smooth functions
$f\in C^{\infty}(\Sigma \times M)$
 which satisfy
     \[
           ||f||_{\epsilon}^2=\sum_{k=0}^{\infty}\epsilon_k
           \inner{\nabla^k f}{\nabla^k f} <\infty,
    \]
where $\nabla^k$ denotes $k$-th hermitian covariant derivative
determined by the metric on the manifold $\Sigma \times M$. This
defines a separable Hilbert space of the subspace of smooth
functions defined on $\Sigma \times M$ and induces topology on the space
$C^{\infty}(\Sigma \times M)$. For a given sequence $\epsilon_k$ we will
call this topology the $\epsilon$-topology, and the space
$C^{\infty}(\Sigma \times M)$ with the $\epsilon$-topology will be
denoted by $V$. Following Floer \cite{F1} one can choose a
sequence $\epsilon_k$ such that the space $V$ is a dense subset of
$L^p(\Sigma \times M)$, for $p>2$.
\begin{proposition}
\label{com}
For every positive number $K$  the $\epsilon$-open set
     \[
           V(K)=\{X \in V\>|\>
           ||f||_{\epsilon}<K \}
     \]
 is relatively compact in the  $C^{\infty}$-topology. Moreover, if
$f_k \to  f_{\infty}$ in $C^{\infty}$-topology and $||f_k||_{\epsilon} \leq K$
for all $k$, then $||f_{\infty}||_{\epsilon} \leq K$.
\end{proposition}
\proof{}
To show that the set $V(K)$ is relatively compact
in the $C^{\infty}$-topology we will use the method of diagonal subsequence.
For each natural $n$ let us introduce a norm $ || . ||_{\epsilon, n} $ by
the formula
     \[
           ||f||_{\epsilon, n}^2=\sum_{k=0}^n \epsilon_k
           \inner{\nabla^k f}{\nabla^k f}.
    \]
Each norm $ || . ||_{\epsilon, n} $ is equivalent to the
corresponding Sobolev norm. Since all functions have support in the
compact set $ \Sigma\times M$ then the natural embedding
      \[
              W(\epsilon, n)\to W(\epsilon, m)
      \]
is compact if $n>m$. Here $W(\epsilon, n)$ denotes the completion in the norm
 $ || . ||_{\epsilon, n} $.

Let $f_k$ be a sequence of smooth functions  such that
$||f_n||_{\epsilon}<K$ for every natural number $k$.
Then it is bounded in the $ || . ||_{\epsilon, 2} $-norm
 and by the above remark
 there exists a  subsequence $f^1_k$ convergent in the
 $ || . ||_{\epsilon, 1} $-norm.
Next, the sequence $f^1_k$ is bounded in
$ || . ||_{\epsilon, 3} $-norm so we can choose
a subsequence $f^2_k$ of the sequence $f^1_k$
convergent in the $ || . ||_{\epsilon, 2} $-norm.

Continuing this process  we will obtain for each $l>1$ a subsequence
 $f^l_k$ of the sequence $f^{l-1}_k$ which is convergent in
 the $ || . ||_{\epsilon, l} $-norm. Choose
the diagonal subsequence $f^k_k$. It is convergent in the
$ || . ||_{\epsilon, l} $-norm for every $l$. Therefore it is convergent in
the $C^{\infty}$-topology. This proves relative compactness in the
$C^{\infty}$-topology.

To show that if $ f_k \to  f_{\infty}$ in $C^{\infty}$-topology and
$||f_k||_{\epsilon} \leq K$ for all $k$, then
$||f_{\infty}||_{\epsilon} \leq K$, it is enough to notice that
      \[
        ||f||_{\epsilon}^2 = \lim_{n \to \infty}
        ||f||_{\epsilon, n} ^2
      \]
for every $f $.
This proves the proposition.

\section{ Universal Moduli Space}
We will study the universal moduli space
  \[
    {\M}(A,J)=\{(u,f)\in {\X}^{1,p}\times   V
    \>\>|\>\>\overline{\partial}_{J,f}(u)=0\},
  \]
where $u$ is solution of the equation (\ref{PDEH}) corresponding
to homology class $A\in H_2(M)$ and to the family of almost
complex structures $J$. Here $V$ denotes the space of functions $f$ 
as described in the Proposition \ref{com}.

Consider the infinite dimensional vector bundle ${\cal E}\to {\X}^{1,p}\times
V$ whose fiber at the point $(u,f)$ is the space
  \[
         {\cal E}(u,f)=L^p(\Omega^{0,1} \otimes_{J} u^*TM)
  \]
of $L^p$-sections of the vector bundle $\Omega^{0,1} \otimes_{J} u^*(TM)$
over $\Sigma$. Then the  moduli space ${\M}(A,J)$
is a zero set ${\cal F}^{-1}(0)$ of the section of the vector bundle
given by the formula
  \[
    {\cal F}(u,f)=\overline{\partial}_{J,f}(u).
  \]
If the point $(u,f)$ is zero of the section ${\cal F}$ then the
differential at this point ,
  \[
    D({\cal F})(u,X):W^{1,p}(u^*TM)\times V
    \to {\cal E}(u,f),
  \]
is given by the formula
 \begin{equation}
 \label{lin}
    D({\cal F})(u,f)(\widehat{u},{\widehat{f}})=D_u\widehat{u}+
    (id,u)^*\left( \overline{\partial}_{\Sigma}(\Phi^{-1}
    \circ \partial_M \widehat{f})\right),
 \end{equation}
since $\overline{\partial}_{\Sigma}$ is linear.

Let $L$ denote the  codimension one  subspace of the space
$V$ orthogonal to the vector $f$ and define the operator
   \[
     D_L({\cal F})(u,f):W^{1,p}(u^*TM)\times L
    \to {\cal E}(u,f)
   \]
as the restriction of the operator $D({\cal F})(u,f)$ to the
subspace $W^{1,p}(u^*TM)\times L$.
 \begin{proposition}
 \label{ontorestricted}
   If the point $(u,f)$ is zero of the section ${\cal F}$
   and $u$ is a  map with $u^*(\omega)\not= 0$ then the linear
   operator $D_L({\cal F})(u,f)$  is onto for suitable chosen family
   of Hermitian connections.
 \end{proposition}
 \proof{}
  We claim that the operators $D({\cal F})(u,f)$ and $ D_L({\cal F})(u,f)$
  have the same range. Indeed, let $\xi= D_u(\widehat{u})+
  \overline{\partial}_{\Sigma}(\Phi^{-1}\circ \partial_M \widehat{f})$. 
 Write $\widehat{f}=tf+g$,
  where $g\in L$. By the Theorem \ref{exact2} we have
   \[
    D_u((id,u)^*\Phi^{-1}\circ \partial_M f)=
   (id,u)^*(\overline{\partial}_{\Sigma}\Phi^{-1}\circ \partial_M f).
   \]
  Therefore,
   \begin{eqnarray*}
    {}\xi&=&D_u(\widehat{u})+(id,u)^*(\overline{\partial}_{\Sigma}
             (t\Phi^{-1}\circ \partial_M f)+
       (id,u)^*(\overline{\partial}_{\Sigma}(\Phi^{-1}\circ \partial_M g)) \\
    {}&=& D_u(\widehat{u}+ (id,u)^*(t\Phi^{-1}\circ \partial_M f))+
         (id,u)^*(\overline{\partial}_{\Sigma}(\Phi^{-1}\circ \partial_M g)) \\
    {}&=&D_L({\cal F})(u,f)(\widehat{u}+ (id,u)^*(t\Phi^{-1}\circ \partial_M f),
        g).
   \end{eqnarray*}
  Thus it is enough to show that the range of $D({\cal F})(u,f)$ is equal
  to $L^p(\Omega^{0,1} \otimes_{J} u^*TM)$.
  Since $D_u$ is a Fredholm operator, by the Remark \ref{fredholm}
  the operator $D({\cal F})$ has a closed range and it is enough to prove
  that the range is dense. Using the Hahn-Banach Theorem it is enough to show
  that  if $\eta \in L^q(\Omega^{0,1} \otimes_{J} u^*TM)$ with
  ${1 \over p}+ {1 \over q} =1 $ satisfies
     \[
       \int \inner{\eta}{D_u \widehat{u}}=0
     \]
  and
     \begin{equation}\label{pt}
       \int \inner{\eta}{\overline{\partial}_{\Sigma}(\Phi^{-1}\circ 
       \partial_M \widehat{f})}=0
     \end{equation}
  for every $\widehat{u} \in W^{1,p}(u^*TM)$ and every
  $\widehat{f}) \in V ,$
  then $\eta \equiv 0$.  From the first equation we obtain that $\eta$
  is a week solution of $D_u^* \eta=0$. However, the coefficients of
  the first order terms of $D_u$ are of class $C^{\infty}$ and the same
  is true for the adjoint $D_u^*$. Thus by elliptic
  regularity $\eta$ satisfies the equation $D_u^* \eta=0$ in the strong
  sense and, moreover, $\eta$ is of class $C^{\infty}$. Hence we can write
     \[
      D_u D_u^* \eta=\Delta \eta+{\mbox{\rm lower order terms}}=0
     \]
  and using the Aronszajn's theorem \cite{Ar} it is enough to show that
  $\eta$ vanishes at some open set.

  By the assumption there is an open set $U\subset \Sigma$ such that the map
  $u$ restricted to the set $U$ is an embedding. Choose $z_0 \in U$
  such that $\eta(z_0) \not= 0$. Let $Y$ be a vector field on $M$ with
  support in a small neighborhood of $u(z_0)$. Choose polar coordinates
  $z=r\exp{(2\pi i \theta)}$ on $\Sigma$ with the property that
  if $Y(r\exp{(2\pi\theta)})\not=0$ then $r_1<r<r_2$ for some positive
  $r_1$ and $r_2$. Next, choose a function $g:(0,\infty)\to \mathbb{R}$ with
  a compact support such that $\int_0^{\infty}g(r)dr=0$, and
  $g(r)>0$ for $r_1<r<r_2$.

  Let $f(r)=\int_0^r g(s)ds$. After complexification it means
  $\overline{\partial}_{\Sigma}(f)=gd\overline{z}$ and
  $\overline{\partial}_{\Sigma}(Yf)= gd\overline{z}\otimes_J Y$.
  We may assume that if $ gd\overline{z}\otimes_J Y(z,u(z))\not=0$
  then
   \[
     \inner{\eta}{ gd\overline{z}\otimes_J Y}(z)>0.
   \]
  Thus if (\ref{pt}) holds $\eta$ must be zero.

 Now we are ready to prove the following theorem
  \begin{theorem}
   \label{wkoncu}
    Let $(u,f)$ with $||f||_{\epsilon}=K$ satisfies
    ${\cal F}(u,X)=0$ and $du$ is of maximum rank at some point.
    Then there is a pair $(u_1,f_1)$
    with $||f_1||_{\epsilon}<K$ such that ${\cal F}(u_1,f_1)=0$.
  \end{theorem}
 \proof{}
 We will need the following version of the implicit function theorem for
 Banach spaces.
   \begin{theorem}
    \label{ift}
     Let $f:E_1\times E_2 \to F$ be a smooth map between Banach spaces.
     Assume that the partial derivative $D_1 f$ is surjective at the point
     $(e_1,e_2) \in E_1\times E_2 $ and admits a bounded right inverse.
     Then for every $f_2\in E_2$ near  $e_2$ there exists $f_1\in E_1$
     such that $f(e_1,e_2)=f(f_1,f_2).$
   \end{theorem}
 Let $L$ denote a tangent space at $X$ to the sphere $S(K)$
 of radius $K$ in the Hilbert space $V$. Then by the Proposition
 \ref{ontorestricted} the partial derivative of ${\cal F}$ at
 $(u,f)$ in the direction $W^{1,p}(u^*TM) \times L$ is onto.
 It has also a bounded right inverse since its restriction to the space
 $W^{1,p}(u^*TM)$ is Fredholm by the Remark \ref{fredholm}.
 Thus we apply the Implicit Function Theorem \ref{ift} to obtain a pair
 $(u_1,f_1)$ with $||f_1||_{\epsilon}<K$ such that ${\cal F}(u_1,f_1)=0$.

  \begin{definition}
    \label{ma}
     The {\bf {extended universal moduli space}} is the space
       \[
          \M_{ex}(A,J)\subset W^{1,p}(u^*TM)\times V \times \mathbb{R}
       \]
      of all triples $(u,f,\lambda)$ with $f\in V$
     such that
      $(u,f)$ is a solution of the equation ${\cal F}(u,f)=0$
    such that homology  class $[u]$ is equal to $A$ and $||f||_{\epsilon}
    =\lambda$.
  \end{definition}
The Theorem \ref{wkoncu} implies
\begin{proposition}
\label{La}
Let $\Lambda (A,J)$ be a set defined as the image of the projection
\[
\pi:\M_{ex}(A,J) \to \mathbb{R}
\]
on the third factor. Then $\Lambda (A,J)$ is open in the set of
all real numbers $\mathbb{R}$.
\end{proposition}

%
\section{Nonnegative Properties of  Symplectic \\
            Area}
%
%
In this section we will prove  the main theorem of these notes. Consider
a solution $u$ of the nonlinear partial differential equation
\begin{equation}
  du+J \circ du \circ j + \overline{\partial}_{\Sigma}(X)(u) =0,
\label{PDEl}
\end{equation}
with properly exact perturbation term.
   \begin{theorem}
     Let $u$ be a solution of the equation (\ref{PDEl}) for
    some $J\in {\cal J}$.
    Then the symplectic form $\omega$ evaluated on the class $[u]$ is
    nonnegative:
     \[
      \int_{\Sigma} u^*\omega \geq 0.
     \]
   \end{theorem}
 \proof{}
 We will prove the theorem by arriving at a contradiction.
 Thus assume that
    \[
      \int_{\Sigma} u^*\omega  <0.
    \]
 Without loss of generality we way assume that $||f||_{\epsilon}=1$.
 Let us  choose a constant $h$
 such that $\omega (C)>h$ for any $J(z)$-holomorphic sphere $C$ in $M$

  Consider the set
    \[
         \Lambda_0=\Lambda(A,J)\cap [0,1].
    \]
 where $A$ is a homology class of $u$, $A=[u]$. (See notation
 of the Proposition \ref{La}).
 Then the set $\Lambda_0$ is not empty since $1\in \Lambda_0$.
 By the Proposition \ref{La} it is also open in $[0,1]$. Let
     \[
        \lambda_1=\inf \Lambda_0.
     \]
 By the definition of $\lambda_1$ there exists a sequence
 $(u_n,f_n,\lambda_n)$ of elements of the extended
 universal moduli space ${\cal M}_{ex}(A,J)$ (Definition \ref{ma})
 such that $\lambda_n \to \lambda_1$.

 By the compactness properties (Proposition \ref{com}) we can assume
 that $f_n \to f^1 $, where $f^1$ satisfies
    \[
       ||f1||_{\epsilon}\leq \lambda_1.
    \]
 By the compactness Theorem
 \ref{Gromov2} there exists a subsequence of
 $(u_n,f_n,\lambda_n,)$ and the collection
 $(u^1; C_1^1,...,C_{r(1)}^1)$
 where $u^1$ is a solution of the equation (\ref{PDEl}) with
 $||f^1||_{\epsilon}\leq \lambda_1$ and $C_1^1,...,C_{r(1)}^1$ are
 $J(z_i)$-holomorphic spheres.
 We have
   \[
       [u]=[u^1]+[C_1^1]+...+[C_{r(1)}^1]
   \]
 so
  \[
   \widetilde{\omega}([\widetilde{u}])=\widetilde{\omega}([{\widetilde{u}}^1])
   +\omega([C_1^1])+...+\omega([C_{r(1)}^1]).
  \]
 See the Remark \ref{tilde}. We have
  \[
   \omega([u^1])\leq \omega([u])<0.
  \]
 For the energy of ${\widetilde{u}}^1$,
 $\widetilde{\omega}([{\widetilde{u}}^1])$, we obtain the following equality:
  \[
   \widetilde{\omega}([{\widetilde{u}}^1])= \widetilde{\omega}(A) -
   \omega([C_1^1])-...-\omega([C_r^1]).
  \]
 We claim that $\lambda_1 > 0$ since otherwise $u^1$ there would be a
 holomorphic curve with negative symplectic area which is impossible.
 Moreover we claim that $r(1)$ ,the number of $J(z_i)$-holomorphic
 spheres, is strictly positive since otherwise
 because of the inequality
  \[
      ||f^1||_{\epsilon}\leq \lambda_1
  \]
 the Theorem \ref{wkoncu} would imply that there was an element of the
 universal moduli space
   \[
    {\M}_{ex}(A,J)
   \]
 corresponding to the parameter $\lambda^1$. But this would contradict the
 fact that $\Lambda_0$ is an open set.

 Assume that for a natural number $k \geq 1 $ we have the following
 data which we will call $k$-data:
  \begin{enumerate}
   \item
    \label{kdataone}
     For each $1\leq i \leq k$ there exists a triple
       \[
       (u^i,f^i,\lambda^i) \in {\M}_{ex}(A^i, J).
       \]
   \item
     For each $1<i \leq k$ numbers $\lambda^i>0$ satisfies
    \begin{enumerate}
     \item
      \[
       \lambda^{i-1} > \lambda^i,
      \]
     \item
      \[
        \lambda^i= \inf \Lambda_{i-1},
      \]
      where
       \[
       \Lambda_{i-1}=\Lambda(A^{i-1},J) \cap [0,1],
      \]
    \end{enumerate}
  \item
 For each $1<i \leq k$ there exists a natural number $r(i) >0$ and a sequence
 $(C_1^i, C_2^i, ... , C_{r(i)}^i )$ of non constant
 $J(z)$-holomorphic spheres such that:
  \begin{enumerate}
    \item
    \label{om}
    \[
      \omega(A^{i-1})=\omega(A^i)+\omega([C_1^i])+...+\omega([C_{r(i)}^i]),
    \]
   \item
    \label{en}
    \[
     \widetilde{\omega}({\widetilde{u}}^i) = \widetilde{\omega}(A^{i-1})
      - \omega([C_1^i])-...-\omega([C_{r(i)}^i]).
    \]
   \end{enumerate}
  \end{enumerate}
 If $\lambda^k >0$ then, using the above method, we can produce
 the $(k+1)$-data.

 Assume that this process never stops i.e. for any natural number $n \in \N$
 there is the $n$-data.

 Then, by the definition of the constant $h$ and by the condition (3)
 of n-data, we obtain
  \begin{eqnarray*}
   {} \widetilde{\omega}({\widetilde{u}}^n)
     & \leq  &  \widetilde{\omega}(A)-\omega(C_1^{1})-...-\omega(C_1^n) \\
   {}               &< &\widetilde{\omega}(A) -nh.
  \end{eqnarray*}
 Choose
  \[
         n> \frac{\widetilde{\omega}(A)}{h}.
  \]
 Then we obtain that $\widetilde{\omega}({\widetilde{u}}^n)<0$. But this
 is a contradiction since the energy of a holomorphic curve can never
 be negative. Therefore, the process must stop somewhere i.e.
  \[
            \lambda^m=0
  \]
 for some $m$.

 However this can not happen by the same reason as ${\lambda}^1$ could
 never be zero and we have arrived at a contradiction.
 Therefore, $\int_{\Sigma} u^*\omega  \geq 0$, and this completes the proof
 of the theorem.


\end{document}